\documentclass[12pt]{article}
\usepackage{amssymb,latexsym}
\usepackage{amsmath,bbm,amsthm}
\usepackage{latexsym}
\usepackage{euscript}
\usepackage[active]{srcltx}

\oddsidemargin=-0.1cm \tolerance 9000 \theoremstyle{definition}

\def\finproof{\hfill\hbox{\vrule width1.0ex height1.5ex}\vspace{2mm}}

\begin{document}

\begin{center}

{{
\textbf{ \sc
Non-local PDEs 
with a state-dependent delay term presented by Stieltjes integral}}}

\vskip7mm

{
\textsc{Alexander V. Rezounenko }}

\smallskip

Department of Mechanics and Mathematics \\ Kharkov University, 4,
Svobody Sqr., Kharkov, 61077, Ukraine \\ 
rezounenko@univer.kharkov.ua

\end{center}




Parabolic partial differential equations with state-dependent delays
(SDDs) are investigated. The delay term presented by Stieltjes
integral simultaneously includes discrete and distributed SDDs. The
singular Lebesgue-Stieltjes measure is also admissible. The
conditions for the corresponding initial value problem to be
well-posed are presented. The existence of a compact global
attractor is proved.

\marginpar{\tiny Nov 27, 2009}

\section{Introduction}\marginpar{sdd11}

We investigate parabolic partial differential  equations (PDEs) with
delay. Studying of this type of equations is based on the
well-developed theories of the ordinary differential equations
(ODEs) with delays \cite{Hale_book,Walther_book,Azbelev} and PDEs
without delays
\cite{Hadamard-1902,Hadamard-1932,Lions,Lions-Magenes-book}. Under
certain assumptions both types of equations  describe a kind of
dynamical systems that are infinite-dimensional, see
\cite{Babin-Vishik,Temam_book,Chueshov_book} and references therein;
see also \cite{travis_webb,Chueshov-JSM-1992,Cras-1995,NA-1998} and
to the monograph \cite{Wu_book} that are close to our work.

In many evolution systems arising in applications the presented
delays are frequently \textit{state-dependent} (SDDs). The theory of
such equations, especially the ODEs, is rapidly developing and many
deep results have been obtained up to now (see e.g.
\cite{Walther2002,Walther_JDE-2003,Walther_JDDE-2007,Krisztin-2003,Mallet-Paret}
and also the survey paper \cite{Hartung-Krisztin-Walther-Wu-2006}
for details and references).

The PDEs with state-dependent delays were first studied in
\cite{Rezounenko-Wu-2006,Hernandez-2006,Rezounenko_JMAA-2007}. 
An alternative approach to the PDEs with discrete SDDs is proposed
in \cite{Rezounenko_NA-2009}. This approach is based on the
so-called \textit{ignoring condition} \cite{Rezounenko_NA-2009}.
Approaches to equations with discrete and distributed SDDs are
 different. Even in the case of ODEs, the discrete SDD essentially complicates the study
 since, in general, the corresponding nonlinearity is not locally Lipschitz
 continuous on open subsets of the space of continuous functions, and familiar results on
 existence, uniqueness, and dependence of solutions on initial data
 and parameters from, say \cite{Hale_book,Walther_book} fail
 (see \cite{Winston-1970} for an example of the non-uniqueness and
 \cite{Hartung-Krisztin-Walther-Wu-2006} for more details).

In this work, in contrast to previous investigations, we consider a
model where two different types of SDDs 
(discrete and distributed) are presented simultaneously (by
Stieltjes integral). The singular Lebesgue-Stieltjes measure is also
admissible. Moreover, all the assumptions on the delay (see A1-A5
below) allow the dynamics when along a solution the number and
values of discrete SDDs may change, the whole discrete and/or
distributed delays may vanish, disappear and appear again. This
property allows us to study models where some subsets of the phase
space are described by equations with purely discrete SDDs, others
by equations with purely distributed SDDs and there are subsets
which need the general (combined) type of the delay. A solution
could be in different subsets at different time moments. This
property particularly means that not only the values of the delays
are state-dependent, but the {\it type} of the delay is {\it
state-dependent} as well. We study mild solutions and their
asymptotic properties (the existence of an attractor is proved). The
results could be applied to the
diffusive Nicholson's blowflies equation with SDDs. 


\section{The model with  state-dependent delay and basic properties}


Consider the following non-local partial differential equation with
a  state-dependent delay term $F$ presented by Stieltjes integral
\begin{equation}\label{sdd11-1}
\frac{\partial }{\partial t}u(t,x)+Au(t,x)+du(t,x) = \big( F(u_t)
\big)(x), \end{equation} with
\begin{equation}\label{sdd11-2}
\big( F(u_t) \big)(x)\equiv \int^0_{-r} \left\{ \int_\Omega b\left(
u(t+\theta,y) \right) f(x-y) dy \right\} \cdot d g(\theta, u_t)
,\quad x\in \Omega,
\end{equation}

\noindent where $A$ is a densely-defined self-adjoint positive
linear operator with domain $D(A)\subset L^2(\Omega )$ and compact
resolvent, which means that $A: D(A)\to L^2(\Omega )$ generates an
analytic semigroup, $\Omega \subset \mathbbm{R}^{n_0}$ is  a smooth
bounded domain, $f: \Omega -\Omega \to R$   is a bounded measurable
function,
$b:\mathbbm{R}\to \mathbbm{R}$ stands for a locally Lipschitz 
map, 
$d \in \mathbbm{R}, d \ge 0$, and the function $g: [-r,0]\times
C([-r,0]; L^2(\Omega)) \to [0,r]\subset \mathbbm{R}_{+}$ denotes
\textit{a state-dependent delay}.
Let $C\equiv C([-r,0];L^2(\Omega))$. Norms defined on $L^2(\Omega)$
and $C$ are denoted by $||\cdot ||$ and $||\cdot ||_C$,
respectively, and $\langle \cdot,\cdot\rangle$ stands for the inner
product in $L^2(\Omega)$. As usual for delay equations, we denote
$u_t\equiv u_t(\theta)\equiv u(t+\theta)$ for $\theta \in [-r,0].$

\medskip

We consider equation \eqref{sdd11-1} with the initial condition
\begin{equation}\label{sdd11-ic}
  u|_{[-r,0]}=\varphi \in C\equiv C([-r,0];L^2(\Omega)). 
\end{equation}

We assume the following.

\smallskip

{\bf A1.)} {\it For any $\varphi\in C $, the function $g:
[-r,0]\times C([-r,0]; L^2(\Omega)) \to \mathbbm{R}$ is of bounded
variation on $[-r,0].$ The variation $V^0_{-r}g$ of $g$ is {\tt
uniformly bounded} i.e.
\begin{equation}\label{sdd11-A1}
\exists M_{Vg}>0\, :\, \forall \varphi\in C \quad \Rightarrow \quad
V^0_{-r}g (\varphi) \le M_{Vg}.
\end{equation}
}%

\smallskip

It is well-known that any Lebesgue-Stieltjes measure (associated with $g$) 
may be split into a sum of three measures: discrete, absolutely
continuous and singular ones. We will denote the corresponding
splitting of $g$ as follows

\begin{equation}\label{sdd11-3}
 g(\theta, \varphi) = g_d(\theta, \varphi) + g_{ac}(\theta, \varphi) +g_s(\theta,
 \varphi),
\end{equation}
where $g_d(\theta, \varphi)$ is a step-function, $g_{ac}(\theta,
\varphi)$ is absolutely continuous and $g_s(\theta, \varphi)$ is
singular continuous functions (see \cite{Kolmogorov-Fomin} for more
details). We will also denote the continuous part by $g_c\equiv
g_{ac}+g_s$.

Now we assume

{\bf A2.)} {\it For any $\theta\in [-r,0],$ the functions $g_{ac}$ and $g_{s}$ are 
continuous with respect to their second coordinates i.e. $\forall
\theta \in [-r,0] \quad \forall \varphi^n,\varphi \in C :
||\varphi^n -\varphi ||_C \to 0\, (n\to +\infty)\Rightarrow
g_{ac}(\theta,\varphi^n) \to g_{ac}(\theta,\varphi)$ and
$g_{s}(\theta,\varphi^n) \to g_{s}(\theta,\varphi).$
}%

\medskip

{\bf Remark~1}. {\it We notice that a discrete state-dependent delay
does not satisfy assumption A2). More precisely, we may consider the
discrete SDD $\eta : C \to [0,r]$ which is presented by the
step-function $g(\theta,\varphi) = 0$ for $\theta\in
[-r,-\eta(\varphi)]$ and $g(\theta,\varphi) = 1$ for $\theta\in
(-\eta(\varphi),0].$ It is easy to see that for any sequence
$\{\varphi^n\}\subset C,$ such that $\eta(\varphi^n) \to
\eta(\varphi)$ and $\eta(\varphi^n) > \eta(\varphi)$ one has for the
value $\theta_0=-\eta(\varphi)$ that $g(\theta_0,\varphi^n)\equiv
1\neq 0\equiv g(\theta_0,\varphi),$ i.e. A2) does not hold.
}%

\medskip

 {\bf A3.)} 
 {\it The step-function $g_d(\theta,
\varphi)$ is continuous with
 respect to its second coordinate in the sense that
 discontinuities of $g_d(\theta, \varphi)$ at points $\{ \theta_k\}\subset
 [-r,0]$ satisfy the property:
  there are {\tt continuous} functions  $\eta_k : C\to [0,r]$ and $h_k : C\to R$ such that $\theta_k=-\eta_k(\varphi)$
 and $h_k (\varphi)$ is the jump of $g_d$ at point
$\theta_k=-\eta_k(\varphi)$ i.e $h_k(\varphi)\equiv g_d(
\theta_k+0,\varphi) - g_d( \theta_k-0,\varphi)$.

Taking into account that $g_d$ may, in general, have infinite number
of points of discontinuity $\{ \theta_k\}$, we assume that the
series $\sum_k  h_k (\varphi) $ converges {\tt absolutely} and {\tt
uniformly on any bounded subsets of C}.
}%

\medskip

{\bf Remark~2}. {\it Following notations of (\ref{sdd11-3}), we
conclude that A3 means that for any $\chi\in C$ one has
$\Phi_d(\chi)\equiv \int^0_{-r} \chi (\theta) \, d
g_d(\theta,\varphi) = \sum_k \chi(\theta_k) \cdot h_k (\varphi)=
\sum_k \chi(-\eta_k(\varphi)) \cdot h_k (\varphi) $.
}%

\medskip

{\bf Lemma~1}. {\it Assume the function $b$ is a Lipschitz map
$(|b(s)-b(t)|\le L_b |s-t|)$, satisfying $ |b(s)|\le C_1|s| +C_2,
\forall s\in \mathbbm{R}$  with $C_i\ge 0$ and $f$ is measurable and
bounded ($|f(x)|\le M_f$). Under assumptions A1)- A3), the nonlinear
mapping $F: C\to L^2(\Omega)$, defined by (\ref{sdd11-2}), is
continuous.}

\medskip

{\bf Remark~3}. {\it We emphasize that nonlinear map $F$ is not
Lipschitz in the presence of discrete SDDs (i.e. when $g\neq g_c$).
The proof of lemma~1 is based on 
  the properties of the uniformly convergent
series and  the first Helly's theorem \cite[page
359]{Kolmogorov-Fomin}.
}

\smallskip

{\it Proof of lemma 1.} We first split our $g$ on continuous
$g_c\equiv g_{ac}+g_s$ and discontinuous $g_d$ parts (see
(\ref{sdd11-3})). This splitting gives the corresponding splitting
of $F=F_c+F_d,$ where $F_c$ corresponds to the continuous part
$g_c\equiv
 g_{ac}+g_s$.


$\lozenge$ Let us first consider the {\it part $F_c$}

We write
\begin{equation}\label{sdd11-4}
F_c(\varphi) - F_c(\psi) = I_1 + I_2,
\end{equation}
where we denote
\begin{equation}\label{sdd11-5}
I_1 =I_1 (x) \equiv \int^0_{-r} \left\{ \int_\Omega \left[ b(\varphi
(\theta,y))-b(\psi(\theta,y))\right] f(x-y)\, dy \right\} \, dg_c
(\theta,\varphi),
\end{equation}

\begin{equation}\label{sdd11-6}
I_2 =I_2 (x) \equiv \int^0_{-r} \left\{ \int_\Omega
b(\psi(\theta,y)) f(x-y)\, dy\right\} \, d\, [
 g_c (\theta,\varphi) - g_c (\theta,\psi)], \quad x\in \Omega.
\end{equation}
One can check that
\begin{equation}\label{sdd11-7}
||I_1|| \le L_b M_f |\Omega|\cdot ||\varphi - \psi ||_C \cdot
V^0_{-r} g (\varphi).
\end{equation}
This estimate and A1 show that $||I_1||\to 0$ when $||\varphi -\psi
||_C \to 0.$ To show that $||I_2||\to 0$ when $||\varphi -\psi ||_C
\to 0$ we use assumptions A1 and A2 to apply the first Helly's
theorem \cite[page 359]{Kolmogorov-Fomin}.

$\lozenge$ Now we prove the continuity of $F_d$ ({\it discrete
delays}).
Let us fix any $\varphi\in C$ and consider a sequence $\{ \varphi^n
\}\subset C$ such that $||\varphi^n-\varphi ||_C \to 0$ when $n\to
\infty$. Our goal is to prove that $||F_d(\varphi^n)-F_d(\varphi) ||
\to 0$.

Following the notations of A3 (see also remark~2), we write
$$F_d(\varphi)=F_d(\varphi)(x)= \sum_k \int_\Omega b(\varphi
(-\eta_k(\varphi), y)) f(x-y) dy\cdot h_k(\varphi)$$ and split as
follows
$$F_d(\varphi^n)-F_d(\varphi) \equiv K^n_1 + K^n_2 +K^n_3,$$
where
$$K^n_1=K^n_1(x)\equiv \sum_k \int_\Omega \left[ b(\varphi^n
(-\eta_k(\varphi^n), y)) - b(\varphi (-\eta_k(\varphi^n), y))\right]
f(x-y) dy\cdot h_k(\varphi^n),$$
$$K^n_2=K^n_2(x)\equiv \sum_k \int_\Omega   b(\varphi (-\eta_k(\varphi^n), y)) f(x-y)
dy\cdot \left[ h_k(\varphi^n)-h_k(\varphi)\right],$$
$$K^n_3=K^n_3(x)\equiv \sum_k \int_\Omega \left[ b(\varphi
(-\eta_k(\varphi^n), y)) - b(\varphi (-\eta_k(\varphi), y))\right]
f(x-y) dy\cdot h_k(\varphi).$$
Using the Lipschitz property of $b$ one may check that
\begin{equation}\label{sdd11-12}
||K^n_1|| \le L_b M_f |\Omega|^{3/2} ||\varphi^n-\varphi ||_C \cdot
\sum_k |h_k(\varphi^n)|.
 \end{equation}
Now we discuss $K^n_2.$ The grough condition of $b$ implies
$|b(\varphi (-\eta_k(\varphi^n), y)) f(x-y)| \le (C_1 |\varphi
(-\eta_k(\varphi^n),y)| + C_2) M_f$. Hence $|\int_\Omega   b(\varphi
(-\eta_k(\varphi^n), y)) f(x-y) dy| \le C_1 M_f \int_\Omega |\varphi
(-\eta_k(\varphi^n),y)| dy +C_2 M_f |\Omega| \le M_f (C_1
|\Omega|^{1/2} ||\varphi||_C + C_2|\Omega|)$. Here we used the
Cauchy-Schwartz inequality for $\int_\Omega |\varphi
(-\eta_k(\varphi^n),y)| dy \le ||\varphi (-\eta_k(\varphi^n))||
\cdot |\Omega|^{1/2} \le ||\varphi||_C\cdot |\Omega|^{1/2}$. One
sees that
$$ |K^n_2(x)| \le M_f (C_1 |\Omega|^{1/2} ||\varphi||_C + C_2|\Omega|)\sum_k |h_k(\varphi^n)-h_k(\varphi)|. $$
Since the right-hand side of the last estimate is independent of
$x$, we get
 \begin{equation}\label{sdd11-13}
 || K^n_2|| \le M_f (C_1 |\Omega|\cdot ||\varphi||_C + C_2|\Omega|^{3/2})\sum_k
 |h_k(\varphi^n)-h_k(\varphi)|.
 \end{equation}
In a similar way we obtain
\begin{equation}\label{sdd11-14}
 || K^n_3|| \le M_f L_b |\Omega| \sum_k
 |h_k(\varphi)|\cdot ||\varphi
(-\eta_k(\varphi^n)) - \varphi (-\eta_k(\varphi)||.
 \end{equation}
Now we should explain why $||K^n_j||\to 0$ as $n\to\infty$ for
$j=1,2,3.$ The first property $||K^n_1||\to 0$ follows from A3 and
(\ref{sdd11-12}). In (\ref{sdd11-13}), the series converges
uniformly with respect to $n$ since the condition
$||\varphi^n-\varphi ||_C \to 0$  implies that $\{ \varphi,
\varphi^n \}$ is a bounded subset of $C$. Assumption A3 guarantees
that each $|h_k(\varphi^n)-h_k(\varphi)|$ is continuous with respect
to $\varphi^n$ and tends to zero when $n\to\infty$. Due to the
uniform convergence we arrive at $||K^n_2||\to 0$. To show that
$||K^n_3||\to 0$ we also mention that each $|h_k(\varphi)|\cdot
||\varphi (-\eta_k(\varphi^n)) - \varphi (-\eta_k(\varphi)||$ (see
(\ref{sdd11-14})) is continuous with respect to $\varphi^n$ and
tends to zero when $n\to\infty$ due to A3 and the continuity of
$\varphi\in C$. The uniform convergence (w.r.t. $\varphi^n$) of the
series in (\ref{sdd11-14}) follows from the estimate
$|h_k(\varphi)|\cdot ||\varphi (-\eta_k(\varphi^n)) - \varphi
(-\eta_k(\varphi)||\le |h_k(\varphi)|\cdot 2 ||\varphi ||_C$
(the right-hand side is independent of $n$!) and the Weierstrass 
dominant (uniform) convergence theorem. We conclude that
$||K^n_3||\to 0$. Since all $||K^n_j||\to 0$ as $n\to\infty$ for
$j=1,2,3$ we proved the property $||F_d(\varphi^n)-F_d(\varphi) ||
\to 0$.
The proof of lemma~1 is complete. \finproof

\section{Mild solutions}


In our study we use the standard

\smallskip

{\bf Definition~1}. 
 {\it A function $u\in C([-r,T]; L^2(\Omega))$ is called a {\tt mild solution}
 on $[-r,T]$ of the initial value problem (\ref{sdd11-1}), (\ref{sdd11-ic}) if it satisfies
 (\ref{sdd11-ic}) and
 \begin{equation}\label{sdd8-3-1}
u(t)=e^{-A t}\varphi(0) + \int^{t}_0 e^{- A (t-s)} \left\{ F(u_s) -
d \cdot u(s)\right\}\, ds, \quad t\in [0,T].
 \end{equation}
}%

\medskip
\medskip

{\bf Theorem~1}. {\it Under assumptions of lemma~1, initial
value problem (\ref{sdd11-1}), (\ref{sdd11-ic})  possesses a 
mild solution for any $\varphi\in C$.}

\smallskip

The existence of a mild solution is a consequence of the continuity
of $F: C\to L^2(\Omega)$, given by lemma~1, which gives us the
possibility to use the standard method based on the Schauder fixed
point theorem (see e.g. \cite[theorem 2.1, p.46]{Wu_book}). The
solution is also global (is defined for all $t\ge -r$), see e.g.
\cite[theorem 2.3, p. 49]{Wu_book}.

\bigskip

To get the uniqueness of mild solutions we need the following
additional assumptions.

\medskip

{\bf A4.)} {\it The total variation of function $g_c\equiv
g_{ac}+g_{s}$ satisfies the {\tt Lipschitz condition}
\begin{equation}\label{sdd11-lip}
 V^0_{-r} [g_{c}(\cdot, \varphi) - g_{c}(\cdot, \psi)] \le L_{Vg_c}
 ||\varphi-\psi||_C.
\end{equation}
}%

\medskip

{\bf A5.)} {\it Discrete generating function $g_{d}$ satisfies the
{\tt uniform  ignoring condition} i.e.

\begin{itemize}
 \item $\exists \eta_{ign}>0$ such that {\tt all} $\eta_k$ and $h_k$ "{\tt ignore}" values of
 $\varphi(\theta)$ for $\theta\in (-\eta_{ign},0]$ i.e. 
 $$\hskip-12mm \exists\, \eta_{ign}>0 : 
 \forall\varphi^1, \varphi^2\in C :
 \forall\theta\in
 [-r,-\eta_{ign}],\,\Rightarrow \varphi^1(\theta)= \varphi^2(\theta)\quad
   \Longrightarrow
   $$
 $$\eta_k (\varphi^1)=\eta_k (\varphi^2),\, h_k (\varphi^1)= h_k \varphi^2). 
 $$
\end{itemize}

}%

\medskip

{\bf Remark~4}. {\it Assumption A5 is the natural generalization to
the case of multiple discrete state-dependent delays of the ignoring
condition introduced in \cite{Rezounenko_NA-2009}. For more details
and examples see \cite{Rezounenko_NA-2009}.
}%

\medskip

{\bf Theorem~2}. {\it Assume the function $b$ is a Lipschitz map
$(|b(s)-b(t)|\le L_b |s-t|)$, satisfying $ |b(s)|\le M_b,\, \forall
s\in R$  
and $f$ is measurable and bounded ($|f(x)|\le M_f$). Under
assumptions A1)- A5), initial value problem (\ref{sdd11-1}), (\ref{sdd11-ic})  possesses a 
unique mild solution for any $\varphi\in C$. The solution is
continuous with respect to initial data i.e. $||\varphi^n -
\varphi||_C\to 0$ implies  $||u^n_t - u_t||_C\to 0$ for any $t\ge
0.$}

\smallskip

{\it Proof of theorem~2.}
The proof is based on 
 the Gronwall lemma, mean value theorem for the Stieltjes integral,
  properties of $g_d$ due to the ignoring condition and the Lebesgue-Fatou
lemma\cite[p.32]{yosida}.

For the simplicity, we first consider a particular case when the
generating function $g=g_c\equiv g_{ac}+g_{s}$ i.e. it does not
contain the discrete delays.

One can check (see (\ref{sdd11-6})) that 
\begin{equation}\label{sdd11-8}
 ||I_2|| \le M_b M_f
|\Omega|^{3\over 2}\cdot  V^0_{-r} [ g_c (\varphi) -g_c (\psi) ].
\end{equation}
Hence (\ref{sdd11-A1}), (\ref{sdd11-4}), (\ref{sdd11-7}),
(\ref{sdd11-8}) and A4 (see (\ref{sdd11-lip})) imply
\begin{equation}\label{sdd11-9}
||F_c(\varphi)-F_c(\psi)|| \le L_{F_c}||\varphi-\psi||_C\quad
\hbox{with } L_{F_c}\equiv M_f|\Omega|  \left( L_b M_{Vg_c} + M_b
|\Omega|^{1\over 2} L_{Vg_c}\right).
\end{equation}
Hence
$$ ||u^1_t-u^2_t||_C \le ||\varphi -\psi||_C + L_{F_c} \cdot \int^t_0
||u^1_s-u^2_s||_C \, ds.
$$
The last estimate (by the Gronwall lemma) implies
$$ ||u^1_t-u^2_t||_C \le
 e^{L_{F_c} t}\cdot ||\varphi -\psi||_C.
$$
That is
\begin{equation}\label{sdd11-CT}
 ||u^1_t-u^2_t||_C \le
C_T\cdot ||\varphi -\psi||_C, \quad \forall t\in [0,T],\qquad \mbox{
with } \, C_T \equiv e^{L_{F_c} T}. 
\end{equation}
We proved the uniqueness of mild solutions and the continuity with
respect to initial data in the case $g=g_c$.

The second particular case $g=g_d$ (the purely discrete delay) and
only one point of discontinuity has been considered in details in
\cite{Rezounenko_NA-2009}. It was proved in
\cite{Rezounenko_NA-2009} that A5  implies the desired result.



Now we consider the general case (both discrete and continuous
delays, including the case of multiple discrete delays). Consider a
sequence $\{ \varphi^n \}\subset C$ such that $||\varphi^n-\varphi
||_C \to 0$ and denote the corresponding mild solutions by
$u^n(t)=u^n(t;\varphi^n)$ and $u(t)=u(t;\varphi)$. Using the
splitting $F=F_d+F_c$, we have, by definition,
$$u^n(t)-u(t) = e^{-At} (\varphi^n(0)-\varphi(0)) +
\int^t_0 e^{-A(t-\tau)} \left\{F_d(u^n_\tau) - F_d(u_\tau)\right\}\,
d\tau $$ $$ + \int^t_0 e^{-A(t-\tau)} \left\{F_c(u^n_\tau) -
F_c(u_\tau)\right\}\, d\tau.
$$
Using (\ref{sdd11-9}), one gets
$$||u^n(t)-u(t)|| = ||\varphi^n(0)-\varphi(0)| +
\int^t_0 ||F_d(u^n_\tau) - F_d(u_\tau)||\, d\tau
+ L_{F_c}\int^t_0 ||u^n_\tau - u_\tau||_C\, d\tau.
$$
Hence
$$
||u^n_t-u_t||_C = ||\varphi^n-\varphi||_C + \int^t_0 ||F_d(u^n_\tau)
- F_d(u_\tau)||\, d\tau + L_{F_c}\int^t_0 ||u^n_\tau - u_\tau||_C\,
d\tau
$$
\begin{equation}\label{sdd11-10}
= G^n(t) + L_{F_c}\int^t_0 ||u^n_\tau - u_\tau||_C\, d\tau,
\end{equation}
where $G^n(t)\equiv ||\varphi^n-\varphi||_C + \int^t_0
||F_d(u^n_\tau) - F_d(u_\tau)||\, d\tau$ is a nondecreasing
function.

Multiply the last estimate by $e^{-L_{F_c} t}$ to get
$$ {d\over dt} \left( e^{-L_{F_c} t} \int^t_0 ||u^n_\tau - u_\tau||_C\, d\tau\right)
\le e^{-L_{F_c} t} G^n(t),
$$
which, after integration from $0$ to $t$, shows that ($G^n(t)$ is
nondecreasing)
$$ e^{-L_{F_c} t} \int^t_0 ||u^n_\tau - u_\tau||_C\, d\tau \le
\int^t_0 e^{-L_{F_c} \tau} G^n(\tau)\, d\tau \le G^n(t)\int^t_0
e^{-L_{F_c} \tau} \, d\tau = G^n(t) \left( 1-e^{-L_{F_c} t}\right)
L^{-1}_{F_c}.
$$
We have
$$ L_{F_c}\int^t_0 ||u^n_\tau - u_\tau||_C\, d\tau \le G^n(t) \left( e^{L_{F_c}
t}-1 \right).
$$
We substitute the last estimate into (\ref{sdd11-10}) to obtain
\begin{equation}\label{sdd11-11}
||u^n_t-u_t||_C \le  G^n(t)\cdot e^{L_{F_c} t}.
\end{equation}
Now our goal is to show that for any fixed $t\in [0, \eta_{ign})$
one has $G^n(t)\to 0$ when $n\to\infty$ (i.e. $||\varphi^n-\varphi
||_C \to 0$).

Let us consider the extension functions
$$
\overline \varphi(s)\equiv \left[\begin{array}{ll}
  \varphi(s) & s \in [-r, 0]; \\
  \varphi(0) & s\in (0, \eta_{ign}) \\
\end{array}
\right. \quad \hbox{ and }\quad  \overline \varphi^n(s)\equiv
\left[\begin{array}{ll}
  \varphi^n(s) & s \in [-r, 0]; \\
  \varphi^n(0) & s\in (0, \eta_{ign}) \\
\end{array}.
\right.
$$
As in \cite{Rezounenko_NA-2009}, the {\it ignoring condition A5}
implies that for all $t\in [0,\eta_{ign})$ we have
$F_d(u_t)=F_d(\overline \varphi_t)$ and $F_d(u^n_t)=F_d(\overline
\varphi^n_t)$. It is easy to see that the convergence $||\varphi^n
-\varphi||_C\to 0$ implies $||\overline \varphi^n_\tau -\overline
\varphi_\tau||_C\to 0$ for any $\tau\in [0, \eta_{ign}).$ Hence the
continuity of $F_d$ implies $||F_d(\overline \varphi^n_\tau)
-F_d(\overline \varphi_\tau)||\to 0$ for any $\tau\in [0,
\eta_{ign}).$ This allows us to use the Lebesgue-Fatou lemma (see
\cite[p.32]{yosida}) for the scalar function $||F_d(\overline
\varphi^n_\tau) -F_d(\overline \varphi_\tau)||$ to conclude that
$G^n(t)\to 0$ when $n\to\infty$ (for any fixed $t\in [0,
\eta_{ign})$). So, we proved the continuity of the mild solutions
with respect to initial functions for all $t\in [0, \eta_{ign})$.
Particularly, it gives the uniqueness of solutions. For bigger time
values we use the chain rule (by the uniqueness) for steps less than
or equal to, say $\eta_{ign}/2$. More precisely, we denote by
$q\equiv \left[ 2t\over \eta_{ign}\right]$ (here $\left[
\cdot\right]$ is the integer part of a real number) and write
$u(t;\varphi)= u(\underbrace{ %
\eta_{ign}/2; u(\eta_{ign}/2; \ldots }_{ q\,\, times };u(t-q\cdot
\eta_{ign}/2; \varphi)))$. The composition of continuous mappings is
continuous. The proof of theorem~2 is complete. \finproof

\bigskip

In the standard way we define an evolution semigroup $S_t: C\to C$
by the rule
$$
S_t \varphi \equiv u_t,
$$
where $u$ is the unique mild solution of (\ref{sdd11-1}),
(\ref{sdd11-ic}).

\medskip

{\bf Remark 5.} {\it The continuity of $S_t$ with respect to time
variable follows from definition~1 (the solution is a continuous
function $u\in C([-r,T]; L^2(\Omega))$).
This and the continuity of $S_t$ with respect to initial function
(see theorem~2) particularly mean that, under assumptions A1)-A5),
the initial value problem (\ref{sdd11-1}), (\ref{sdd11-ic}) is {\tt
well-posed} in the space $C$ in the sense of J.~Hadamard
\cite{Hadamard-1902,Hadamard-1932}.}

\medskip

The last remark means that the pair $( S_t,C)$ forms the dynamical
system (for the definition see e.g
\cite{Babin-Vishik,Temam_book,Chueshov_book}).

Following the line of argument given in \cite[theorem
2]{Rezounenko_NA-2009} we show that the dynamical system $( S_t,C)$
generated by initial value problem (\ref{sdd11-1}), (\ref{sdd11-ic})
possesses a compact global attractor (for more details on attractors
see, for example \cite{Babin-Vishik,Temam_book,Chueshov_book})).

More precisely, we have the following result.

\medskip

{\bf Theorem~3.} {\it Assume the function $b : \mathbbm{R}\to
\mathbbm{R}$ is a
Lipschitz and bounded map ($|b(w)|\le C_b$ for all $w\in
\mathbbm{R}$) and $f: \Omega -\Omega \to \mathbbm{R}$ is a bounded
and measurable function ($|f(\cdot )|\le M_f$). Let assumptions
A1-A5 be satisfied.
 Then the dynamical
system $( S_t,C )$ has a compact global attractor which is a compact
set in all spaces $C_\delta\equiv C([-r,0]; D(A^\delta)),
\forall\delta\in [0,{1\over 2}).$}

The proof is based on the classical theorem on the existence of a
compact global attractor for a dissipative and asymptotically
compact semigroup~\cite{Babin-Vishik,Temam_book,Chueshov_book} and
technique developed in \cite[theorem 2]{Rezounenko_NA-2009}.

\medskip

As an application we can consider the diffusive Nicholson's
blowflies equation (see e.g. \cite{So-Yang}) with state-dependent
delays, i.e. equation (\ref{sdd11-1}) where $-A$ is the Laplace
operator with the Dirichlet boundary conditions, $\Omega\subset
\mathbbm{R}^{n_0}$ is a bounded domain with a smooth boundary, the
nonlinear (birth) function $b$ is given by $b(w)=p\cdot we^{-w}$.
The function $b$ is bounded, so under assumptions A1-A5, we conclude
that the initial value problem (\ref{sdd11-1}) and (\ref{sdd11-ic})
is well-posed in $C$ and the dynamical system $(S_t,C)$ has a
compact global attractor (theorem~3).



\medskip





\end{document}